\documentclass[a4paper,12pt]{article}
\usepackage{amscd}
\usepackage{amsmath,amsfonts,amssymb,amscd}
\usepackage{indentfirst,graphicx,epsfig}
\usepackage{graphicx,psfrag}
\usepackage{ifpdf}
\usepackage{float}

\input{epsf}

\newtheorem{thm}{Theorem}

\begin{document}
\title{\Large\bf Sharp upper bound for the rainbow connection number of a graph with
diameter 2\footnote{Supported by NSFC No.11071130.}}
\author{\small Jiuying Dong, Xueliang Li\\
\small Center for Combinatorics and LPMC-TJKLC\\
\small Nankai University, Tianjin 300071, China\\
\small jiuyingdong@126.com, lxl@nankai.edu.cn}
\date{}
\maketitle

\begin{abstract}

Let $G$ be a connected graph. The \emph{rainbow connection number
$rc(G)$} of a graph $G$ was recently introduced by Chartrand et al.
Li et al. proved that for every bridgeless graph $G$ with diameter
2, $rc(G)\leq 5$. They gave examples for which $rc(G)\leq 4$.
However, they could not show that the upper bound 5 is sharp.
It is known that for a graph $G$ with diameter 2, to determine
$rc(G)$ is NP-hard. So, it is interesting to know the
best upper bound of $rc(G)$ for such a graph $G$. In
this paper, we use different way to obtain the same upper bound, and
moreover, examples are given to show that the upper is best possible.\\[3mm]


\end{abstract}

\section{Introduction}

All graphs considered in this paper are simple, finite and
undirected. Undefined terminology and notations can be found in
\cite{Bondy-Murty}. Let $G$ be a graph, and $c: E(G) \rightarrow
\{1,2,\cdots,k\}, k\in N$ be an edge-coloring, where adjacent edges
may be colored the same. A graph $G$ is \emph{rainbow connected} if
for every pair of distinct vertices $u$ and $v$ of $G$, $G$ has a
$u-v$ path $P$ whose edges are colored with distinct colors. The
minimum number of colors required to make $G$ rainbow connected is
called the \emph{rainbow connection number} of $G$, denoted by
$rc(G)$. These concepts were introduced in \cite{Chartrand-Johns},
In \cite{Basavaraju-Chandran}, Basavaraju et al. showed that for
every bridgeless graph $G$ with radius $r$, $rc(G)\leq r(r+2)$, and
the bound is sharp. In \cite{S-Ch}, Chakraborty et al. investigated
the hardness and algorithms for the rainbow connection number, and
showed that given a graph $G$, deciding if $rc(G) = 2$ is
NP-Complete. In particular, computing $rc(G)$ is NP-Hard. It is
well-known that almost all graphs have diameter 2. So, it is
interesting to know the best upper bound of $rc(G)$ for a graph $G$
with diameter 2. In \cite{H-Li}, Li et al. proved that for every
bridgeless graph $G$ with diameter 2, $rc(G)\leq 5$. They gave
examples for which $rc(G)\leq 4$. However, they could not show that
the upper bound 5 is sharp. In this paper, we use different way to
obtain the same upper bound, and moreover, examples are given to
show that the upper is best possible. Our main result is the
following:

\begin{thm}\label{thm}
If $G$ is a connected bridgeless graph with diameter 2, then
$rc(G)\leq 5$. Moreover, the upper bound is sharp.
\end{thm}

Before proceeding, we need some notation and terminology. For two
subsets $X$ and $Y$ of $V$, an $(X,Y)$-path is a path which connects
a vertex of $X$ and a vertex of $Y$, and whose internal vertices
belong to neither $X$ nor $Y$. We use $E[X,Y]$ to denote the set of
edges of $G$ with one end in $X$ and the other end in $Y$, and
$e(X,Y)=|E[X,Y]|$. The $k$-step open neighborhood of $X$ is defined
as $N^{2}(X)=\{v\in V(G)|d(v,X)=2\}$. Let $N[S]=N(S)\cup S$. For a
connected graph $G$, the eccentricity of a vertex $v$ is
$ecc(v)=\max_{x\in V(G)}d_G(v,x)$. The radius of $G$ is
$rad(G)=\min_{x\in V(G)}ecc(x)$. The diameter of $G$ is $\max_{x\in
V(G)}ecc(x)$, denoted by $diam(G)$.

First of all, we give the following example to show that the upper
bound $rc(G)\leq 5$ is sharp.

\section{An example}

At first, we construct a graph $G$ as follows:

Let $P_1=uv_1w_1, P_2=uv_2w_2,\cdots, P_k=uv_kw_k$ be $k$ internally
disjoint paths of length 2 with $k\geq 17$, and when $ i\neq j$,
$w_i\neq w_j, 1 \leq i,j\leq k $. For any two different vertices
$w_i, w_j$, we join $w_i$ to $ w_j$. Thus we get a graph $G$. Let
$V=\{v_1,v_2,\cdots,v_k\}$, $W=\{w_1,w_2,\cdots,w_k\}$. We know that
$V$ is an independent set, $G[W]$ is a complete subgraph of $G$. So
the diameter of $G$ is 2.

In any edge-coloring $c: E(G)\rightarrow \{1,2,3,4\}$ of $G$, each
2-length $u-w_i$ path can be colored in at most 16 different ways.
By the Pigeonhole Principle, there exist $i,j, i\neq j, 1 \leq
i,j\leq k$ such that $c(uv_i)=c(uv_j), c(v_iw_i)=c(v_jw_j)$.
Consider any rainbow path $R$ between $v_i$ and $v_j$. We know that
$R$ can contain only two edges of $\{uv_i,uv_j,v_iw_i,v_jw_j\}$, and
$R$ must pass through another path $P_l,l\neq i, l\neq j$ and one
edge of $\{w_iw_l, w_jw_l\}$. That is, $R=v_iw_iw_lv_luv_j$ or
$R=v_jw_jw_lv_luv_i$. Hence, $rc(G)\geq 5$. Now, we color the edges
of $G$ as follows: Let $c(uv_1)=1, c(v_1w_1)=2, c(uv_i)=3,
c(v_iw_i)=4, 2\leq i\leq k $, for any edge $e\in E(G[W])$, $c(e)=5$.
It is easy to check that for any two vertices $v_1,v_2$ of $G$,
there is a rainbow path connecting them. Hence, we get that
$rc(G)\leq 5$. That is, $rc(G)=5$. Thus, $G$ is a required sharp
example.

\section{Proof of Theorem \ref{thm}}

Let $G$ be a bridgeless graph with $diam(G)=2$, $r$ be the radius of
$G$, $u$ be the center vertex of $G$. We say the colorings of the
following cycles containing $u$ to be $appropriate$: let
$C_3=uv_1v_2u$ be a 3-cycle where $v_1,v_2\in N^{1}(u)$, and let
$c(uv_1)=1, c(uv_2)=2, c(v_1v_2)=3(4)$; let $C_4=uv_1v_2v_3u$ be a
4-cycle where $v_1,v_3\in N^{1}(u), v_2\in N^{2}(u)$, and let
$c(uv_1)=1, c(uv_3)=2, c(v_1v_2)=3,c(v_3v_2)=4$; let
$C_5=uv_1v_2v_3v_4u$ be a 5-cycle where $v_1,v_4\in N^{1}(u),
v_2,v_3\in N^{2}(u)$, and let $c(uv_1)=1, c(uv_4)=2, c(v_1v_2)=3,
c(v_3v_4)=4, c(v_2v_3)=5$. We know that the shortest cycles passing
through $u$ are only the above mentioned $C_3, C_4$ or $C_5$.

Let $B_1,B_2,\cdots, B_b, B_{b+1}, B_{b+2}\subset N^{1}(u) $ satisfy
the following two conditions:\\
(1). For $1\leq i\neq j \leq b+2$, $\bigcup_{i=1}^{b+2} B_i=N^{1}(u)$£¬
$B_i\cap B_j=\emptyset$, and for $1\leq i\leq b,B_i=N_{N^{1}(u)}[b_i]$;\\
(2). If $B_{b+1}\neq\emptyset$, then $B_{b+1}$ is an independent
set, and for any $w\in B_{b+1}, wb_i\not\in E(G)$, but $\exists
w'\in\bigcup_{i=1}^{b} (B_i\setminus\{b_i\})$ such that $ww'\in
E(G)$; if $B_{b+2}\neq\emptyset$, that is to say
$N^{2}(u)\neq\emptyset$, then $B_{b+2}$ is an independent set, and
$E( B_{b+2},X\setminus B_{b+2} )= \emptyset$.

Let $c(ub_i)=1$. For any $e \in E(u,B_i\setminus \{b_i\})$, let
$c(e)=2$. For any $e \in E(u, B_{b+1})$, let $c(e)=1$. In the
following we will construct a series of sets $S_1\subset S_2\subset
S_3\subset\cdots\subset S_k$, where $S_1= \bigcup_{i=1}^{b+1} B_i$,
$N^{1}(u)\subset S_k$, and every vertex $v\in S_i\setminus S_1$ is
in a cycle with the appropriate coloring. Given an $S_i$, if
$N^{1}(u)\subset S_i$, then we stop the procedure by setting $k=i$.
Otherwise, we construct $S_{i+1}\supset S_i$ as follows:

For any $v\in N^{1}(u)\setminus S_i$, we select a cycle $R_i$ such
that $R_i$ is a shortest cycle containing $uv$. \ \ \ \ \ \ \ \ (*)

Subject to (*), we further choose $R_i$ such that $R_i$ contains as
many vertices of $G\setminus S_i$ as possible.

If $V(R_i)\cap S_i=\{u\}$, then we give $R_i$ an $appropriate$
coloring. Otherwise, $V(R_i)\cap (S_i\setminus\{u\})\neq\emptyset$.
There will appear two cases:

(I). $R_i$ is a 4-cycle. Let $R_i=uv_1v_2v_3u$. Either $v_2\not\in
S_i, v_3\in S_i$, or $v_2\in S_i,v_3\in S_i$. In the former case, we
color the edges of $R_i$ as follows: if $c(uv_3)=2$, then let $
c(v_3v_2)=4, c(v_2v_1)=3, c(v_1u)=1$; if $c(uv_3)=1$, then let $
c(v_3v_2)=3, c(v_2v_1)=4, c(v_1u)=2$. In the latter case, we color
the edges of $R_i$ as follows: if $c(uv_3)=1,c(v_3v_2)=3$, then let
$c(uv_1)=2, c(v_1v_2)=4$; if $c(uv_3)=2,c(v_3v_2)=4$, then let
$c(uv)=1, c(v_1v_2)=3$.

(II). $R_i$ is a 5-cycle. Let $R_i=uv_1v_2v_3v_4u$. If $v_2\in S_i$,
then $\exists v'\in S_i\cap N^{1}(u)$ such that $uv_1v_2vu$ is a
4-cycle, a contradiction to $|R_i|=5$. So $v_2\not\in S_i$. Either
$v_4\in S_i, v_2, v_3\not\in S_i$, or $v_3,v_4\in S_i$. In the
former case, we color the edges of $R_i$ as follows: if $c(uv_4)=1$,
then let $c(v_4v_3)=3, c(v_2v_3)=5, c(v_1v_2)=4, c(uv_1)=2$; if
$c(uv_4)=2$, then let $c(v_4v_3)=4, c(v_2v_3)=5, c(v_1v_2)=3,
c(uv_1)=1$. In the latter case, we color the edges of $R_i$ as
follows: if $c(uv_4)=2, c(v_4v_3)=4$, then let $c(v_2v_3)=5,
c(v_1v_2)=3, c(uv_1)=1$; if $c(uv_4)=1,c(v_4v_3)=3$, then let
$c(v_2v_3)=5, c(v_1v_2)=4, c(uv_1)=2$.

When the above procedure ends, we have constructed $S_k$, and
$N^{1}(u)\subset S_k$. If $N^{2}(u)\subset S_k$, then for any $e\in
E(G[N^{1}(u)])$, we let $c(e)=3$, and we will color the remaining
uncolored edges by a used color. From the above construction, we can
see that any vertex of $S_k$ is in a shortest cycle with the
appropriate coloring. Hence for any two vertices $v,v'\in S_k$, if
$v,v'$ are in the same cycle with the appropriate coloring, then it
is obvious that there is a rainbow path connecting $v,v'$. If $v$ is
in some cycle $R$ with the appropriate coloring, $v'$ is in another
cycle $R'$ with the appropriate coloring, because $R, R'$ both
contain $u$, it is not difficult to check that there is a rainbow
path connecting $v,v'$.

Otherwise, $N^{2}(u)\setminus S_k\neq \emptyset$. Let
$N^{1}(u)=X\cup Y$, where for any $x\in X$, we let $c(ux)=1$, for
any $y\in Y$ we let $c(uy)=2$. Let all $S, T, Q \subseteq
N^{2}(u)$ be as large as possible which satisfy the following conditions:\\
(1). For any $s\in S, E(s,X)\neq \emptyset$ but $E(s,Y)= \emptyset$;
For any $t\in T, E(t,Y)\neq \emptyset$ but $E(t,X)=
\emptyset$; For any $q\in Q, E(q,X)\neq \emptyset$ and $E(q,Y)\neq \emptyset$.\\
(2). For any $s\in S, E(s,T\cup Q)\neq\emptyset$; For any $t\in T,
E(t,S\cup Q)\neq\emptyset$.

For any $e\in E(S,X)$, let $c(e)=3$, for any $e\in E(T,Y)$, let
$c(e)=4$, for any $e\in E(Q,X)$, let $c(e)=3$, and for any $e\in
E(Q,Y)$, let $c(e)=4$. If $N^{2}(u)= S\cup T\cup Q$, then for any $
e \in E(G[N^{1}(u)])$, let $c(e)=3$, for any $e\in E(S, T\cup Q)$,
let $c(e)=5$, and we use a used color to color the remaining
uncolored edges. We may see that any vertex of $S_k\cup S\cup T\cup
Q$ is in some cycle with the appropriate coloring. Hence for any two
vertices of $S_k\cup S\cup T\cup Q$, there is a rainbow path
connecting them. That is to say, $G$ is rainbow connected.
Otherwise, let $ N^{2}(u)\setminus (S\cup T\cup Q)=P\cup L$, where
for any $ p\in P, E(p,X)\neq\emptyset$, for any $\ell\in L,
E(\ell,Y)\neq\emptyset$. By the maximality of $S,T,Q$, we know that
$E(P,Y)=\emptyset, E(P, T\cup Q)=\emptyset, E(L,X)=\emptyset, E(L,
S)=\emptyset$. Suppose that $P \neq\emptyset$ and $L\neq\emptyset$,
then $E(P,L)=\emptyset$, and for any $p\in P, \ell\in L$, the
distance of $p$ and $\ell$ is more than 2, a contradiction to
$Diam(G)=2$. Hence, without loss of generality, in the following we
may assume $P \neq\emptyset, L=\emptyset$.

Let $P=P_1\cup P_2$ where for any $p_1\in P_1, e(p_1,X)=1$, for any
$p_2\in P_2, e(p_2,X)\geq 2$. So we may get that for any $y\in Y,
x_{p_1}y\in E(G)$ where $x_{p_1}\in X, x_{p_1}p_1\in E(G)$, and for
some given $y\in Y,\exists x_{p_2}\in X$ such that $p_2x_{p_2}\in
E(G), x_{p_2}y\in E(G)$. Hence, if $P\neq\emptyset$, then by the
construction of $S_k$, we may get that for any $v\in B_{b+2},
c(uv)=1$. Otherwise, if $\exists v\in B_{b+2}$ such that $c(uv)=2$,
then for any $p\in P, \exists x_p\in X $ such that $px_p,x_pv\in
E(G)$, thus $v\in N^{1}(u)\setminus B_{b+2}$, a contradiction.
Hence, we may get $S_1\neq\emptyset$. If $P_1\neq\emptyset$, then
$E(P_1,B_{b+2})=\emptyset$. Otherwise, if $\exists p_1\in P_1,
x_{b+2}\in B_{b+2}$ such that $p_1x_{b+2}\in E(G)$, then for any
$y\in Y,yx_{b+2}\in E(G)$, thus $ x_{b+2}\in N^{1}(u) \setminus
B_{b+2}$, a contradiction. If $P_2\neq\emptyset$, then for any $p\in
P_2, e(p,B_{b+2})\leq 1$. If $\exists p_2\in P_2$ such that
$e(p,B_{b+2})\geq 2$, then by the construction and the coloring of
$S_k$, $\exists x_{b+2},x_{b+2}'\in B_{b+2}$ such that $u,
x_{b+2},x_{b+2}'$ are in a 4-cycle with the appropriate coloring,
and we may assume $c(ux_{b+2})=1, c(ux_{b+2}')=2$, a contradiction
to the fact that for any $x\in B_{b+2}, c(ux)=1$. In the following
we divide into two cases to show that $G$ is rainbow connected.

{\bf Case 1.} $P_1=\emptyset$.

For any $e\in E(X,Y)$, let $c(e)=3$, for any $e\in E(P, S)$, let
$c(e)=2$, and for any $e\in E(P,P)$, let $c(e)=5$. Let $y_p\in Y$,
for any $p\in P_2,\exists x_p\in X\setminus B_{b+2} $ such that
$px_p,x_py_p\in E(G)$, let $c(px_p)=5$. We use color 4 to color the
remaining uncolored edges of $E(p,X)$, and use a used color to color
the remaining uncolored edges of $E(G)$.  In the following we show
that $G$ is rainbow connected. For $(p,p')\in (P_2\times P_2)$,
either $px_py_pux_{p'}p'$ is a rainbow path, where $x_{p'}\in X,
c(p'x_p')=4$, or $px_pp'$ is a rainbow path, where $ c(x_pp')=4$.
For each pair $(p,t)\in (P_2\times T), px_puy_tt$ is a rainbow path,
where $y_t\in Y$. For each pair $(p,q)\in (P_2\times Q), px_puy_qq$
is a rainbow path, where $y_q\in Y$, (Latter, we will omit these
notes). For each pair $(p,y)\in (P_2\times Y)$, $px_puy $ is a
rainbow path. For each pair $(p,x)\in (P_2\times X)$, $px_py_pux $
is a rainbow path. For each pair $(p,s)\in (P_2\times S)$, either
$ps$, or $px_ss$, or $pp's$ is a rainbow path. For each pair
$(x,x')\in (X\setminus B_{b+2}\times X), xy_xux'$ is a rainbow path.
By the construction of $S_k$, for any $x\in  B_{b+2}$, either $x$ is
in a 4-cycle with the appropriate coloring, or $x$ is in a 5-cycle
with the appropriate coloring. Hence for each pair $(x,x')\in
(B_{b+2}\times X)$, either $xqy_qux'$ is a rainbow path, or
$xsv_sy_sux'$ is a rainbow path. Therefore for any two vertices of
$G$, there is a rainbow path connecting them.

{\bf Case 2.} $ P_1\neq\emptyset$.

{\bf Subcase 2.1.}  $|X|=1$.

Let $D_1,D_2,\cdots, D_d,D_{d+1}\subset P$ satisfy the following conditions:\\
(1). For $1\leq i\neq j \leq d+1,D_i\cap D_j=\emptyset$, and for
$1\leq i\leq d, D_i=N_P[d_i]$.\\
(2). If $D_{d+1}\neq \emptyset$, then $D_{d+1}$ is an independent
set, and for any $v\in D_{d+1}, vd_i\not\in E(G)$, and $ \exists
v'\in N^{2}(u)\setminus\bigcup_{i=1}^{d} d_i $ such that $vv'\in
E(G)$.

We let $c(xd_i)=1$, for any $e \in E(x, D_i\setminus\{d_i\})$ let
$c(e)=2$, for any $e \in E(G[P])$ let say $c(e)=3$, for any $e\in
E(x, D_{d+1})$ let $ c(e)=1$, for any $e\in E(D_{d+1},
N^{2}(u)\setminus D_{d+1}) $ let $ c(e)=4$, and for any $e\in
E(G[N(u)])$ let $c(e)=3$. Then we use a used color to color the
remaining uncolored edges. For any two vertices of $P\setminus
D_{d+1}$, it is not difficult to check that there is a rainbow path
connecting them. For each pair $(v_1, v_2)\in (D_{d+1}\times
D_{d+1})$, either $\exists v\in \bigcup_{i=1}^{d} D_i\setminus
\{d_i\}$ such that $v_1vxv_2$ is a rainbow path, or $\exists s\in S$
such that $v_1sxv_2$ is a rainbow path. For any $p\in P$, all $pxy,
pxq, pxst$ are rainbow paths. Hence $G$ is rainbow connected.

{\bf Subcase 2.2.} $|X|\geq 2$.

{\bf Subsubcase 2.2.1.} $\exists x\in X\setminus B_{b+2} $ such that
$E(x,P_1)=\emptyset$.

Let $X_1\subseteq X$ be the set of all vertices which are adjacent
to the vertices of $P_1$. Let $X_2=X\setminus (X_1\cup B_{b+2})$,
$x_2\in X_2$, and let $y_{x_2}\in Y$ such that $x_2y_{x_2}\in E(G)$.
Let $P_1'\subset P $ be the set of all vertices which are adjacent
to the vertices of $X_1$, and let $P_2'=P\setminus P_1'$. For any
$e\in E(X_1, Y)$ let $c(e)=3$, for any $e\in E(X_2, Y)$ let
$c(e)=4$, for any $e\in E(X_1, P_1')$ let $c(e)=1$, and for any
$e\in E(x_2, P)$ let $c(e)=5$. For any $p\in P_2'$, if $px_2\in
E(G)$, we use color 3 to color the remaining uncolored edges of
$E(p,X)$. If $px_2\not\in E(G)$, we know that $\exists x_p\in X_2$
such that $p x_p\in E(G)$, and let $c(p x_p)=5$. Then we use color 3
to color the remaining uncolored edges of $E(p,X)$. For any $e\in
E(G[P])$ let $c(e)=2$, and for any $e\in E(P, S)$ let $c(e)=2$.
Finally, we use a used color to color the remaining uncolored edges
of $E(G)$. For each pair $(p,p')\in (P_1'\times P_1')$, if $px_2\in
E(G)$, then $px_2y_{x_2}x_{p'}p'$ is a rainbow path; if $p'x_2\in
E(G)$, then $p'x_2y_{x_2}x_{p}p$ is a rainbow path; if $px_2\not\in
E(G), p'x_2\not\in E(G)$, then either $\exists p''\in P$ such that
$pp''x_2y_{x_2}x_{p'}p'$ is a rainbow path, or $\exists s\in S$ such
that $psv_sy_sx_{p'}p'$ is a rainbow path, where $v_s\in T\cup Q,
y_s\in Y$. For each pair $(p_1,p_2)\in (P_1'\times P_2')$,
$p_1x_{p_1}y_{p_2}x_{p_2}p_2$ is a rainbow path, where
$c(x_{p_2}p_2)=5$.  For each pair $(p,x)\in (P\times X)$, either
$px$, or $\exists p'\in P$ such that $pp'x$ is a rainbow path, or
$\exists s\in S$ such that $psx$ is a rainbow path. For any $y\in Y,
s\in S, t\in T, q\in Q$, all $p_1x_{p_1}y, p_1x_{p_1}y_sv_ss,
p_1x_{p_1}y_tt, p_1x_{p_1}y_qq $ are rainbow paths. For each pair
$(p,p')\in (P_2'\times P_2')$, $ px_py_pux_p'p'$ is a rainbow path,
where $c(px_p)=5, c(p'x_p')=3$. For any $ y\in Y, s\in S, t\in T,
q\in Q $, all $px_py_pux, px_puy, px_py_pux_ss, px_puy_tt, px_puy_qq
$ are rainbow paths. For each pair $(x,x')\in (X\times X)$, similar
to Case 1, there is a rainbow path connecting them. Therefore, we
conclude that $G$ is rainbow connected.

{\bf Subsubcase 2.2.2.}  For any $ x\in X, E(x,P_1)\neq\emptyset$.

Let $x_1\in X, y\in Y$. For any $E(x_1, P)$ let $c(e)=5$, for any
$E(X, P)\setminus E(x_1, P)$ let $c(e)=3$, for any $E(x_1, Y)$ let
$c(e)=4$, for any $E(X\setminus \{x_1\}, Y)$ let $c(e)=1$, for any
$e\in E(G[P])$ let $c(e)=2$, and for any $e\in E(P, S)$ let
$c(e)=2$. Finally, we use a used color to color the remaining
uncolored edges of $E(G)$. For each pair $(p,p')\in (P\times P)$, if
$px_1\not\in E(G), p'x_1\not\in E(G)$, then either $\exists p''\in
P$ such that $pp''x_1yx_{p'}p'$ is a rainbow path, or $\exists s\in
S$ such that $ps,sx_1\in E(G)$, thus $psv_sy_sx_{p'}p'$ is a rainbow
path; if $px_1\in E(G), p'x_1\not\in E(G)$, then $px_1yx_{p'}p'$ is
a rainbow path; if $p'x_1\in E(G), px_1\not\in E(G)$, then
$p'x_1yx_{p}p$ is a rainbow path; if $px_1\in E(G), p'x_1\in E(G)$,
then for any $ x\in X\setminus ( B_{b+2}\cup\{x_1\} )$, either
$pxyx_1p'$ is a rainbow path, or $\exists p''\in P$ such that
$pp''xyx_{1}p'$ is a rainbow path, or $\exists s\in S$ such that
$psxyx_{1}p'$ is a rainbow path. For each pair $(p,s)\in (P\times
S)$, either $px_1yux_{s}s$ is a rainbow path, or $px_py_sv_ss$ is a
rainbow path. For each pair $(p,t)\in (P\times T)$, either
$px_1uy_tt$ is a rainbow path, or $px_puy_tt$ is a rainbow path. For
each pair $(p,q)\in (P\times Q)$, either $px_1uy_qq$ is a rainbow
path, or $px_puy_qq$ is a rainbow path. For each pair $(p,x)\in
(P\times X)$, either $px$, or $\exists p'\in P$ such that $pp'x$ is
a rainbow path, or $\exists s\in S$ such that $psx$ is a rainbow
path. For each pair $(p,y)\in (P\times Y)$, either $ px_1y$, or
$px_py$ is a rainbow path. For any $x\in X\setminus\{x_1\},x_1yux$
is a rainbow path. For any $x,x'\in X\setminus
(B_{b+2}\cup\{x_1\})$, let $p\in P$ such that $px\in E(G)$. Then
either $xpx_1yx'$ is a rainbow path, or $\exists p'\in P$ such that
$xpp'x_1yx'$ is a rainbow path, or $\exists s\in S$ such that
$xpsv_sy_sx'$ is a rainbow path. For each pair $(x,x')\in
(B_{b+2}\times X)$, similar to Case 1, there is a rainbow path
connecting them. Therefore, we conclude that $G$ is rainbow
connected. This completes the proof of Theorem 1.

\end{document}